\documentclass[a4paper,12pt]{article}
\usepackage[dvips]{epsfig}
\usepackage{amsmath,amssymb, amsbsy, amstext,amscd,amsfonts}
\usepackage{color,enumerate,euscript,graphicx,hyperref}
\usepackage{srcltx,tikz}
\input amssym.def
\input amssym.tex

\newtheorem{lem}{Lemma}[section]%
\newtheorem{theorem}[lem]{Theorem}%
\newtheorem{exam}[lem]{Example}%
\newtheorem{prop}[lem]{Proposition}%
\newtheorem{conj}[lem]{Conjecture}%

\def\a{\alpha}

 \def\O{\Omega}

 \def\og{\overline G} \def\oh{\overline H}

 \def\ox{\overline X}  \def\o1{\overline 1}

\def\olg{\overline g} \def\ola{\overline a} 
  \def\olc{\overline c} \def\olb{\overline b}

\def\o{\overline}   \def\olb{\overline b}
\def\di{\bigm|} \def\lg{\langle} \def\rg{\rangle}

\def\Aut{\hbox{\rm Aut\,}}  \def\Syl{\hbox{\rm Syl}}
  
  \def\mod{\hbox{\rm mod }}

 \def\AGL{\hbox{\rm AGL}} \def\GL{\hbox{\rm GL}}  \def\P\GL{\hbox{\rm P\GL}}
\def\GF{\hbox{\rm GF}}  \def\FF{{\hbox{\sf F\kern-.43emF}}}
 \def\char{\hbox{\rm char}}

\def\o{\hbox{\rm o}}

\def\char{ \, {\rm char}\,}

\def\ZZ{\mathbb{Z}}  

\def\nd{\mathrel{\bigm|\kern-.7em/}} 
 \def\f{\noindent}
\def\qed{\hfill $\Box$} \def\demo{\f {\bf Proof}\hskip10pt}

\begin{document}
\begin{center}
{\bf\large  The Product of a Semi Dihedral Group And a Cyclic Group}
\end{center}


\begin{center}
Hao Yu{\small \footnotemark}\\
\medskip
 {\small
Capital Normal University,\\ School of Mathematical Sciences,\\
Beijing 100048, People's Republic of China
}
\end{center}

\footnotetext{Corresponding author: 3485676673@qq.com. }

\renewcommand{\thefootnote}{\empty}
\begin{abstract}
  Let $X(G)=GC$ be a group, where $G$ is a semi dihedral group and $C$ is a cyclic group such that
  $G\cap C=1$. In this paper, $X(G)$ will be characterized.
  \end{abstract}

\section{Introduction}
A group $G$ is said to be properly {\it factorizable} if $G=AB$ for two proper subgroups $A$ and $B$ of $G$, while the  expression $G=AB$ is called a {\it factorization} of $G$.
Furthermore, if $A\cap B=1$,  then  we say that $G$ has an {\it exact factorization}.

Factorizations of groups naturally arise from the well-known  Frattini's argument, including  its version in permutation groups.
One of the most famous results about factorized groups might be one of theorems of It\^o, saying that  any group is metabelian whenever it is the product of two abelian subgroups (see \cite{Ito1}).
Later, Wielandt and  Kegel showed that the product of two nilpotent subgroups must be soluble
(see \cite{Wie1958} and \cite{Ke1961}).
Douglas showed that  the product of two cyclic groups must be super-solvable (see \cite{D1961}).
The factorizations of the finite almost simple groups were determined in \cite{LPS1990}  and
the  factorizations of almost simple groups with a solvable factor  were determined in  \cite{LX2022}.
There are many other papers related to factorizations, for instance, finite products of soluble groups, factorizations with one nilpotent factor and so on.
Here we are not able to list all references and the readers may refer to a survey paper \cite{AK2009}.

In this paper, we shall focus on the product group $X=X(G)=GC$,
for a finite group $G$ and a cyclic group $C$ such that $G\cap C=1$.
Recently, $X(G)$, where $G$ is either a dihedral group or a generalized quaternion group,
has been characterized (\cite{DYL2023}).
The core-free property of $C$ is frequently referenced in this article.
Suppose that $C$ is core-free. Then $X$ is  also called a {\it skew product group} of $G$.
Recall that the skew morphism of a group $G$ and a skew product group $X$ of $G$
were introduced by Jajcay and \v{S}ir\'a\v{n} in \cite{JS2002},
which is related to the studies of regular Cayley maps of $G$.
For the reason of the length of the paper, we are not able to explain them in detail.
Recently, there have been a lot of results on skew product groups $X$ of some particular groups $G$.
(1) Cyclic groups: So far there exists no classification of such product groups.
For  partial results, see~\cite{CJT2016,CT,DH,KN1,KN2,Kwo}.
(2) Elementary abelian $p$-groups: a global structure  was characterized in \cite{DYL}.
(3) Finite nonabelian simple group or finite nonabelian characteristically simple groups:
they were classified  in  \cite{BCV2019} and \cite{CDL}, respectively.
(4) Dihedral groups: Based on big efforts of several authors working on regular Cayley maps (see \cite{CJT2016,HKK2022,KKF2006, KMM2013,KK2017,KK2016,RSJTW2005,KK2021,WF2005,WHY2019,Zhang2015,
Zhang20152,ZD2016}),
the final classification of  skew product groups of dihedral groups was given in \cite{HKK2022}.
(5) Generalized quaternion groups: they were classified  in \cite{HR2022}, \cite{KO2008} and \cite{DYL2023}.

Let $X(G)=GC$ be a group, where $G$ is a semi dihedral group and $C$ is a cyclic group such that $G\cap C=1$.
In this paper, we shall give a characterization for $X(G)$ and some property of $X(G)$.

Throughout this paper, set
\begin{eqnarray}\label{main0}
\begin{array}{ll}
&C=\lg c\di c^m=1\rg\cong\ZZ_m, \,m\geq2,\\
&D=\lg a, b\di a^{n}=b^2=1,\,a^b=a^{-1}\rg\cong D_{2n},\,n\ge 2.\\
&Q=\lg a, b\di a^{2n}=1,\,b^2=a^n,a^b=a^{-1}\rg\cong Q_{4n},\,n\ge 2.\\
&SD=\lg a, b\di a^{4n}=b^2=1,\,a^b=a^{2n-1}\rg\cong SD_{8n},\,n\ge 2.
\end{array}
\end{eqnarray}
Let $G=SD$ and $X=X(G)=GC=\lg a, b\rg \lg c\rg$.
Then $\lg a\rg \lg c\rg $ is unnecessarily a subgroup of $X$.
Clearly, $X$ contains a subgroup $M$ of the biggest order such that
$\lg c\rg\le M\subseteqq\lg a\rg\lg c\rg$.
This subgroup $M$ will play an important role in this paper.
From now on by $S_X$ we denote the core $\cap_{x\in X}S^x$ of $X$ in a subgroup $S$ of $X$.

There are two main theorems in this manuscript.
In Theorem~\ref{main1}, the global structure of our group $X$ is characterized.
\begin{theorem} \label{main1}
Let $G=SD$ and $X=X(G)=G\lg c\rg$, where $\o(c)=m\ge 2$ and $G\cap \lg c\rg=1$.
Let $M$ be the subgroup of the biggest order in $X$ such that
$\lg c\rg \le M\subseteqq \lg a\rg \lg c\rg$.  Then one of items in Tables \ref{table1} holds.
\begin{table}
  \center \caption {The forms of $M$, $M_X$ and $X/M_X$}\label{table1}
  \begin{tabular}{cccc}
  \hline
  Case &$M$ & $M_X$  & $X/M_X$\\
  \hline
   1 & $\lg a\rg\lg c\rg$   & $\lg a\rg\lg c\rg$     &   $\ZZ_2$ \\
   2 &$\lg a^2\rg \lg c\rg$ & $\lg a^2\rg\lg c^2\rg$ &   $D_8$   \\
   3 &$\lg a^2\rg \lg c\rg$ & $\lg a^2\rg\lg c^3\rg$ &   $A_4$   \\
   4 &$\lg a^4\rg \lg c\rg$ & $\lg a^4\rg\lg c^3\rg$ &   $S_4$   \\
   5 &$\lg a^4\rg \lg c\rg$ & $\lg a^4\rg\lg c^4\rg$ &   $(\ZZ_4\times\ZZ_2)\cdot\ZZ_4$   \\
   6 &$\lg a^3\rg \lg c\rg$ & $\lg a^3\rg\lg c^4\rg$ &   $S_4$   \\
   \hline
  \end{tabular}
  \end{table}
\end{theorem}

Clearly, $M$ is a product of two cyclic subgroups,
 which  has not been determined so far, as mentioned before,
However, further  properties of our group $X$ is given in Theorem~\ref{main2}.

\begin{theorem} \label{main2}
Let $G=SD$ and $X=X(G)$, and $M$ defined as above.
Then we have $\lg a^4, c\rg \le C_X(\lg c\rg_X)$, $|X: C_X(\lg c\rg_X)|\le 8$ if $n$ is odd,
$\lg a^2, c\rg \le C_X(\lg c\rg_X)$, $|X: C_X(\lg c\rg_X)|\le 4$ if $n$ is even.
Moreover, if $\lg c\rg _X=1$, then $M_X\cap \lg a^4\rg \lhd M_X$.
In particular, if $\lg c\rg _X=1$ and $M=\lg a\rg \lg c\rg$, then $\lg a^4\rg \lhd X$.
\end{theorem}

After this introductory section, some  preliminary results will be given in Section 2, Theorems~\ref{main1} and \ref{main2} will be proved in Sections 3 and 4, respectively.

\section{Preliminaries}
In this section, the notation and elementary facts used in this paper are  collected.

\subsection{Notation}
In this paper, all the groups are supposed to be finite.
We set up the notation below, where $G$ and $H$ are groups, $M$ is a subgroup of $G$,
$n$ is a positive integer and $p$ is a prime number.
\begin{enumerate}
  \setlength{\itemsep}{0ex}
  \setlength{\itemindent}{-0.5em}
  \item[] $|G|$ and $\o(g)$: the order of $G$ and an element $g$ in $G$, resp.;
  \item[] $H\leq G$ and $H<G$: $H$ is a subgroup of $G$ and $H$ is a proper subgroup of $G$, resp.;
  \item[] $[G:H]$: the set of cosets of   $G$ relative to a subgroup $H$;
  \item[] $H\lhd G$ and $H\char~G$: $H$ is a normal and characteristic subgroup of $G$, resp.;
  \item[] $G'$ and $Z(G)$: the derived subgroup and the center of $G$ resp.;
  \item[] $M_G:=\cap_{g\in G} M^g$, the core of $M$ in $G$;
  \item[] $G\rtimes H$: a semidirect product of $G$ by $H$, in which $G$ is  normal;
  \item[] $G.H$:  an extension of $G$ by $H$, where $G$ is normal;
  \item[] $C_M(G)$: centralizer of $M$ in $G$;
  \item[] $N_M(G)$: normalizer of $M$ in $G$;
  \item[] $\Syl_p(G)$: the set of all Sylow $p$-subgroups of $G$;
  \item[] $[a,b]:=a^{-1}b^{-1}ab$, the commutator of $a$ and $b$ in $G$;
  \item[] $\O_1(G)$: the subgroup $\lg g\in G\di g^p=1\rg$ of $G$ where $G$ is a $p$-group;
  \item[] $\mho_n(G)$: the subgroup $\lg g^{p^n}\di g\in G$ of $G$ where $G$ is a $p$-group;
  \item[] $S_n$: the symmetric group of degree $n$ (naturally acting on $\{1, 2, \cdots, n\}$);
  \item[] $A_n$: the alternating group of degree $n$ (naturally acting on $\{1, 2, \cdots, n\}$);
  \item[] $\GF(q)$: finite field of $q$ elements;
  \item[] $\AGL(n,p)$: the affine group on $\GF^n(q)$.
\end{enumerate}
\subsection{Elementary facts}
\begin{prop}\cite[Theorem 1]{M1974}\label{solvable}
The finite group $G=AB$ is solvable,
where both $A$ and $B$ are subgroups with cyclic subgroups of index no more than 2.
\end{prop}

Recall that a group $H$ is said a {\it Burnside group}
if every permutation group containing a regular subgroup isomorphic to $H$ is
either 2-transitive or imprimitive.
The following results are well-known.

\begin{prop} \cite[Theorem 25.3 and Theorem 25.6]{W1964} \label{Burnside}
Every cyclic group of composite order is a Brunside group.
Every dihedral group is a Burnside group.
\end{prop}

\begin{prop}\cite[Corollary 1.2]{Y2023}\label{SDB}
Every semi dihedral group of order $8n\, (n\geq3)$ is a Brunside group.
\end{prop}

\begin{prop}\label{complement}\cite[Satz 1]{G1952}
Let $N\leq M\leq G$ such that $(|N|,|G:M|)=1$ and $N$ be  an abelian normal subgroup of $G$.
If $N$ has a complement in $M$, then $N$ also has a complement in $G$.
\end{prop}

\begin{prop}\cite[Theorem 4.5]{H1967}\label{NC}
Let $H$ be the subgroup of $G$. Then $N_G(H)/C_G(H)$ is isomorphic to a subgroup of $\Aut(H)$.
\end{prop}

\begin{prop}\cite[Theorem]{Luc} \label{cyclic}
If $G$ is a transitive permutation group of degree $n$  with a cyclic point-stabilizer,
then $|G|\le n(n-1)$.
\end{prop}

\begin{prop} \cite[Satz 1 and Satz 2]{Ito1} \label{mateabel}
Let $G=AB$ be a group, where both $A$ and $B$ are abelian subgroups of $G$. Then
\begin{enumerate}
  \item[\rm(1)]  $G$ is meta-abelian, that is, $G'$ is abelian;
  \item[\rm(2)]  if $G\ne 1$, then  $A$ or $B$ contains   a normal subgroup $N\ne 1$ of $G$.
\end{enumerate}
\end{prop}

\begin{prop}\cite[Theorem 11.5]{H1967}\label{ab}\label{matecyclic}
Let $G=\lg a\rg\lg b\rg$ be a group.
If $|\lg a\rg|\leq |\lg b\rg|$, then $\lg b\rg_G\ne1$.
If both $\lg a\rg$ and $\lg b\rg$ are $p$-groups where $p$ is an odd prime, then $G$ is matecyclic.
\end{prop}

\begin{prop}\cite[Corollary 1.3.3]{W1999}\label{sylowp}
Let $G=AB$ be a group, where both $A$ and $B$ are subgroups of $G$.
And let $A_p$ and $B_p$ be Sylow $p$-subgroups of $A$ and $B$ separately, for some prime $p$.
Then $A_pB_p$ is the Sylow $p$-subgroup of $G$.
\end{prop}

\begin{prop}\cite[Theorem 12.5.1]{H1959}\label{pcyclic}
Let $p$ is an odd prime.
Then every finite $p$-group $G$ containing  a cyclic maximal subgroup is isomorphic to
(1)  $\ZZ_{p^n}$;\,   (2) $\lg a,b\di a^{p^{n-1}}=b^p=1,\,[a,b]=1\rg,\,n\geq2$; or
(3) $\lg a,b\di a^{p^{n-1}}=b^p=1,\,[a,b]=a^{p^{n-2}}\rg,\,n\geq3$.
\end{prop}

\begin{prop}\cite[Lemma 4.1]{DMM2008}\label{AGL}
Let $n\ge 2$ be an integer and $p$ a prime.
Then $\AGL(n, p)$ contains an element of order $p^n$ if and only if $(n, p)=(2, 2)$ and $\AGL(2, 2)\cong S_4$.
\end{prop}

Recall that our group $X(D)=DC$ and $X(Q)=QC$, where $D$ is a dihedral group of order $2n$,
$Q$ is a generalized quaternion group of order $4n$ and
$C$ is a cyclic group of order $m$ such that $D\cap C=1$, where $n, m\ge 2$.
Then we have the following results.

\begin{lem}\label{cd}
 Suppose that $X(D)$ is a solvable and has a faithful 2-transitive permutation  representation relative to a subgroup $M$, which is of index  a composite order.
 Then $X(D)\leq\AGL(k,p)$. Moreover,
\begin{enumerate}
  \item[\rm(i)]  if $X(D)$ contains an element of order $p^k$, then $X(D)=S_4$;
  \item[\rm(ii)] if the hypotheses holds for $M=C$ where $C$ is core-free, then $X(D)=A_4$.
\end{enumerate}
\end{lem}

\begin{lem} \label{DM}
Let $G\in\{ Q,D\}$ and $X=X(G)=\lg a, b\rg\lg c\rg$, and
let $M$ be the subgroup of the biggest order in $X$ such that $\lg c\rg\le M\subseteqq\lg a\rg\lg c\rg$.  Then one of items in Tables \ref{tableD1} holds.
\begin{table}
  \center   \caption {The forms of $M$, $M_X$ and $X/M_X$}\label{tableD1}
  \begin{tabular}{cccc}
  \hline
  Case &$M$ & $M_X$  & $X/M_X$\\
  \hline
   1 & $\lg a\rg\lg c\rg$   & $\lg a\rg\lg c\rg$     &   $\ZZ_2$ \\
   2 &$\lg a^2\rg \lg c\rg$ & $\lg a^2\rg\lg c^2\rg$ &   $D_8$   \\
   3 &$\lg a^2\rg \lg c\rg$ & $\lg a^2\rg\lg c^3\rg$ &   $A_4$   \\
   4 &$\lg a^4\rg \lg c\rg$ & $\lg a^4\rg\lg c^3\rg$ &   $S_4$   \\
   5 &$\lg a^3\rg \lg c\rg$ & $\lg a^3\rg\lg c^4\rg$ &   $S_4$   \\
   \hline
  \end{tabular}
  \end{table}
\end{lem}
\begin{lem}\label{Da^2}
Suppose that $G=D,\,X=X(D),\,\lg a\rg \lg c\rg\leq X$ and $\lg c\rg_X=1$.
Then $\lg a^2\rg\lhd X$.
\end{lem}
\vskip 3mm
Recall that our group $X(SD)=SDC$, where $SD$ is a dihedral group of order $8n$ and
$C$ is a cyclic group of order $m$ such that $SD\cap C=1$, where $n, m\ge 2$.
Then we have the following result.

\begin{lem}\label{cd1}
 Suppose that $X(SD)$ is a solvable and has a faithful 2-transitive permutation  representation relative to a subgroup $M$, which is of index  a composite order.
 Then $X(SD)\leq\AGL(k,p)$. Moreover,
\begin{enumerate}
  \item[\rm(i)]  if $X(SD)$ contains no element of order $p^k$;
  \item[\rm(ii)] if the hypotheses holds for $M=C$ where $C$ is core-free, then there no exists $X(SD)$.
\end{enumerate}
\end{lem}
\demo Set $\O=[X(SD):M]$.  Let $N$ be a minimal normal subgroup of $X(SD)$.
Since $X(SD)$ is solvable, $N\cong \ZZ_p^k$ for some prime $p$ and integer $k$.
Since $X(SD)$ is 2-transitive, it is  primitive, which implies that $N$ is transitive on $\O$ and so is regular on $\O$.
Therefore, $X(SD)=N\rtimes X(D)_\a \le \AGL(k,p)$, for some $\a\in \O$.
Since $X(SD)$ is 2-transitive and $|\O|=p^k$, we know $|X(D)_\a|\ge p^k-1$ for any $\a\in \O$.

(i) Arguing by contradiction, assume that $X(SD)$ contains an element of order $p^k$. By Proportion~\ref{AGL}, we  get $(k,p)=(2,2)$ so that  $X(SD)=S_4$, reminding $|\O|$ is not a prime.
But it is impossible, a contradiction.

(ii) Arguing by contradiction, assume that $X(SD)$ exists.
Let $M=C$ where $C$ is core-free.
Set $C=\lg c\rg $ and $o(c)=m$.
Then $X(SD)=N\rtimes \lg c\rg $, where $\lg c\rg $ is a Singer subgroup of $\GL(k,p)$.
Then  both $SD$ and $N$ are regular subgroups, that is  $|SD|=8n=|\O|=p^k$, which implies $p=2$.
Now, we have $|X(SD)|=2^{k}(2^k-1)=8n\cdot m=p^k\cdot m$ and so $m=2^k-1$.
Since both $N$ and $SD$ are Sylow 2-subgroups of $X(SD)$ and $N\lhd X(SD)$, we get $SD=N$.
Since $SD$ is not an abelian group, we get a contradiction.
\qed
\vskip 3mm
Using Lemmas~\ref{SDB} and \ref{cd1}, we get the following lemma.
\begin{lem}\label{SDI}
Every semi dihedral group of order $8n\, (n\geq3)$ and
let $X$ be a permutation group containing a regular subgroup isomorphic to $G$.
If $X$ contains a point stable subgroup of a cycle, then $X$ is   imprimitive.
\end{lem}
\vskip 3mm
To show Theorems~\ref{main1} and \ref{main2}, we need some examples.
\begin{exam}\label{SD16or24}
Let $G=SD_{16}$ or $SD_{24}$ and $X(G)$.
If $\lg c\rg$ is core-free and $G\ntriangleleft X(G)$, then in the isomorphic sense, we have the following result.

\f(i) Suppose that $G=SD_{16}$.

\begin{tabular}{cccc}
  \hline
  Case &$X(G)$\\
  \hline
   1 &$\lg a,b,c\di a^8=b^2=c^3=1,\,a^b=a^3,\,a^c=a^7c^2,\,b^c=ac^2\rg$ \\
   2 &$\lg a,b,c\di a^8=b^2=c^4=1,\,a^b=a^3,\,a^c=a^3,\,b^c=abc^2\rg$ \\
   3 &$\lg a,b,c\di a^8=b^2=c^6=1,\,a^b=a^3,\,a^c=bc^5,\,b^c=a^3c^5\rg$ \\
   4 &$\lg a,b,c\di a^8=b^2=c^6=1,\,a^b=a^3,\,a^c=a^2bc^2,\,b^c=a^4bc^2\rg$ \\
   5 &$\lg a,b,c\di a^8=b^2=c^8=1,\,a^b=a^3,\,a^c=a^5,\,b^c=abc^4\rg$ \\
   6 &$\lg a,b,c\di a^8=b^2=c^8=1,\,a^b=a^3,\,a^c=a,\,b^c=abc^4\rg$ \\
   7 &$\lg a,b,c\di a^8=b^2=c^8=1,\,a^b=a^3,\,a^c=abc^2,\,b^c=a^2bc^4\rg$ \\
   8 &$\lg a,b,c\di a^8=b^2=c^8=1,\,a^b=a^3,\,a^c=abc^2,\,b^c=a^6bc^4\rg$ \\
   9 &$\lg a,b,c\di a^8=b^2=c^8=1,\,a^b=a^3,\,a^c=a^2bc^2,\,b^c=ac^6\rg$ \\
   10 &$\lg a,b,c\di a^8=b^2=c^8=1,\,a^b=a^3,\,a^c=a^2bc^2,\,b^c=a^5c^6\rg$ \\
  \hline
\end{tabular}

\f(ii) Suppose that $G=SD_{24}$.

\begin{tabular}{cccc}
  \hline
  Case &$X(G)$\\
  \hline
   1 &$\lg a,b,c\di a^{12}=b^2=c^3=1,\,a^b=a^5,\,a^c=a^9c^2,\,b^c=bc^2\rg$ \\
   2 &$\lg a,b,c\di a^{12}=b^2=c^3=1,\,a^b=a^5,\,a^c=a^9c^2,\,b^c=a^8b\rg$ \\
   3 &$\lg a,b,c\di a^{12}=b^2=c^4=1,\,a^b=a^5,\,a^c=bc^3,\,b^c=a^8b\rg$ \\
   4 &$\lg a,b,c\di a^{12}=b^2=c^4=1,\,a^b=a^5,\,a^c=a,\,b^c=a^3bc^2\rg$ \\
   5 &$\lg a,b,c\di a^{12}=b^2=c^4=1,\,a^b=a^5,\,a^c=a^5,\,b^c=abc^2\rg$ \\
   6 &$\lg a,b,c\di a^{12}=b^2=c^6=1,\,a^b=a^5,\,a^c=abc^2,\,b^c=a^8bc^4\rg$ \\
   7 &$\lg a,b,c\di a^{12}=b^2=c^6=1,\,a^b=a^5,\,a^c=a^3bc^5,\,b^c=a^4c^3\rg$ \\
   8 &$\lg a,b,c\di a^{12}=b^2=c^6=1,\,a^b=a^5,\,a^c=a^9c^2,\,b^c=a^2b\rg$ \\
   9 &$\lg a,b,c\di a^{12}=b^2=c^6=1,\,a^b=a^5,\,a^c=a^3c^2,\,b^c=bc^2\rg$ \\
   10 &$\lg a,b,c\di a^{12}=b^2=c^6=1,\,a^b=a^5,\,a^c=a^3c^2,\,b^c=a^2b\rg$ \\
   11 &$\lg a,b,c\di a^{12}=b^2=c^{12}=1,\,a^b=a^5,\,a^c=a,\,b^c=abc^6\rg$ \\
   12 &$\lg a,b,c\di a^{12}=b^2=c^{12}=1,\,a^b=a^5,\,a^c=abc^4,\,b^c=a^{10}bc^4\rg$ \\
   13 &$\lg a,b,c\di a^{12}=b^2=c^{12}=1,\,a^b=a^5,\,a^c=bc^{10},\,b^c=a^3c^6\rg$ \\
  \hline
\end{tabular}
\end{exam}

\begin{exam}\label{Z4Z2}
Let $G=\lg a,b\di a^4=b^2=1,a^b=a\rg\cong\ZZ_4\times\ZZ_2$ and $X(G)$.
If $\lg c\rg$ is core-free and $G\ntriangleleft X(G)$, then in the isomorphic sense, we have the following result.
\begin{tabular}{cccc}
  \hline
  Case &$X(G)$\\
  \hline
   1 &$\lg a,b,c\di a^4=b^2=c^2=1,\,a^b=a,\,a^c=a,\,b^c=a^2b\rg$ \\
   2 &$\lg a,b,c\di a^4=b^2=c^2=1,\,a^b=a,\,a^c=a^3,\,b^c=b\rg$ \\
   3 &$\lg a,b,c\di a^4=b^2=c^2=1,\,a^b=a,\,a^c=ab,\,b^c=b\rg$ \\
   4 &$\lg a,b,c\di a^4=b^2=c^4=1,\,a^b=a,\,a^c=a,\,b^c=a^3bc^2\rg$ \\
   5 &$\lg a,b,c\di a^4=b^2=c^4=1,\,a^b=a,\,a^c=abc^2,\,b^c=a^2\rg$ \\
  \hline
\end{tabular}
\end{exam}

\section{Proof of Theorem~\ref{main1}}
\vskip 3mm
To prove Theorem~\ref{main1}, let $G=SD$, defined in Eq(\ref{main0}).
Let $X=X(G)=G\lg c\rg$.
Let $M$ be the subgroup of the biggest order in $X$ such that $\lg c\rg \le M\subseteqq \lg a\rg \lg c\rg$,
and set $M_X=\cap_{x\in X}M^x$.
By Proposition~\ref{solvable}, $X$ is solvable.
\vskip 3mm
Before showing Theorem~\ref{main1}, we shall show a special case which plays an important role in the proof of Theorem~\ref{main1}.
\begin{lem}\label{n=p}
Let $G=SD$, $\lg c\rg\cong\ZZ_m$ and $X=X(G)$.
Suppose that $\lg c\rg_X=1$, $H=\lg a^p,b\rg\lg c\rg\leq X$ and
there exists no nontrivial element $a^j\in X$ such that $\lg a^j\rg\lg c\rg\le X$.
Then $n\neq p$ where $p$ is a prime and more than 3.
\end{lem}
\demo
By Example~\ref{SD16or24}, we know $n>3.$
Arguing by contradiction, assume that $n=p$ is a prime and more than 3.
If $G_X\neq1$, we get that $\lg a\rg_X\neq1$ as $p$ is a prime and more than 3,
which implies $\lg a\rg_X\lg c\rg\le X$, a contradiction.
So in what follows, we assume that $G_X=1$.

Set $\lg c_1\rg:=\lg c\rg_H$ and $a_1:=a^p$.
Consider $\oh=H/\lg c_1\rg=\lg\ola^p,\olb\rg\lg \olc\rg$.
Since $\lg \ola^p,\olb\rg\cong\ZZ_4\times\ZZ_2$ and there exists no nontrivial element $\ola^j\in \oh$ such that $\lg \ola^j\rg\lg \olc\rg\le \oh$, by Example~\ref{Z4Z2}, we get that $\oh\cong(\ZZ_4\times\ZZ_2)\cdot\ZZ_4$ and the forms of $\oh$ only have the following one case:
$$\oh=\lg\ola_1,\olb,\olc\di\ola_1^4=\olb^2=\olc^4=1,\,\ola_1^{\olb}=\ola_1,\,
{\ola_1}^{\olc}=\ola_1\olb\olc^2,\,\olb^{\olc}=\ola_1^2\rg.$$
Then we get $4\di m$.
Let $a_0$ and $c_0$  be involutions in $\lg a\rg$ and $\lg c\rg$, respectively.
Noting that in the above forms, we get
\begin{eqnarray}\label{f1.1}
\olb^{\olc}=\ola_0,\,\ola_0^{\olc}=\olb,\,(\ola_0\olb)^{\olc}=\ola_0\olb.
\end{eqnarray}
Since $\lg c_1\rg\lhd H$, we know $c\in C_H(\lg c_1\rg)\lhd H$.
Then by the definition of $\oh$, we get $\lg a^{2p},b,c\rg\leq C_H(\lg c_1\rg)$.
Then in the perimage $X$, Eq(\ref{f1.1}) corresponds to
\begin{eqnarray}\label{f1.2}
b^c=a_0c_0^i,\, a_0^c=bc_0^j,\,(a_0b)^c=a_0bc_0^k,
\end{eqnarray}
where $i,j,k\in\{ 0,1\}$.
If $c_0\notin\lg c_1\rg$, then we get that $i=j=k$ contradicting with $G_X\neq1$.
So $c_0\in\lg c_1\rg$.
Since $G_X=1$, we know that $(a_0b)^c=a_0bc_0$ and the pair $(i,j)$ is either $(0,1)$ or $(1,0)$.
Then Eq(\ref{f1.2}) is either
$$b^c=a_0,\, a_0^c=bc_0,\,(a_0b)^c=a_0bc_0,$$
or
$$b^c=a_0c_0,\, a_0^c=b,\,(a_0b)^c=a_0bc_0.$$
Then one can check that for any subgroup $N$ of $H$,
if $N\lhd H$ and $a_0\in N$, then $b\in N$.

Consider $H_X$.
Since $a_1\in \cap_{l_3} H^{a^{l_3}}=\cap_{l_1, l_2, l_3} H^{c^{l_1}b^{l_2}a^{l_3}}= H_X$,
we get that $\lg a_1\rg \leq H_X$.
Since $a_0\in H_X$ and $H_X\lhd H$, we get that $b\in H_X$.
Since $H_X\lhd X$, we get $a^2\in H_X$, a contradiction.
\qed
\vskip 3mm
\f{\bf Proof of Theorem~\ref{main1}}
Let $G=SD$ so that $X=X(SD)$.
Remind that $m, n\ge 2$, $|X|$ is even and more than 32.
Let $X$ be a minimal count-example.
Then we shall carry out the proof by the following two steps.
\vskip 3mm
{\it Step 1: Show that $M_X=1$.}
\vskip 3mm
Arguing by contradiction, assume that $M_X\ne 1$.
Set $M=\lg a^i\rg\lg c\rg$ for some $i$.
Since $a^i\in \cap_{l_2, l_3} M^{a^{l_2}b^{l_3}}=\cap_{l_1, l_2, l_3} M^{c^{l_1}a^{l_2}b^{l_3}}= M_X$,
we get that $M_X=M_X\cap (\lg a^i\rg\lg c\rg )=\lg a^i\rg \lg c^r\rg$ for some $r$.
Set $\ox:=X/{M_X}=\og\lg\olc\rg$.
Then we claim that $\og\cap\lg\olc\rg=1$.
In fact, for any $\olg=\olc'\in \og\cap \lg\olc\rg$ for some $g\in G$ and $c'\in \lg\olc\rg$,
we have $gc'^{-1}\in M_X$,
that is $g\in \lg a^i\rg$ and $c'\in \lg c^r\rg$, which implies $\olg=\olc'=1$.
Therefore, $\og\cap\lg\olc\rg=1$.
Let $M_0/M_X=\lg \ola^j\rg \lg \olc\rg $ be the largest subgroup of $\ox$
containing $\lg \olc\rg$ and contained in the subset $\lg \ola\rg \lg \olc\rg $.
Then $\lg \ola^j\rg \lg \olc\rg =\lg \olc\rg \lg \ola^j\rg$.
Since
$$\lg a^j\rg \lg c\rg M_X=\lg a^j\rg M_X\lg c\rg=\lg a^j\rg \lg a^i\rg \lg c\rg \quad {\rm and}\quad \lg c\rg \lg a^j\rg M_X=\lg c\rg M_X\lg a^j\rg =\lg c\rg \lg a^i\rg \lg a^j\rg ,$$
we get $\lg a^i, a^j\rg \lg c\rg \le X$.
By the maximality of $M$, we have $\lg a^i, a^j\rg =\lg a^i\rg $ so  that $M_0=M$.

Consider $\ox$. Note that $\og$ is one of a dihedral group, a generalized quaternion group and a semi dihedral group.
Suppose that $\og$ is either a dihedral group or a generalized quaternion group.
Then noting $M_0/M_X=M/M_X$ is core-free in $\ox$, by Lemma~\ref{DM},
we get that $\ox$ is isomorphic to  $\ZZ_2,\,D_8,\,A_4$ or $S_4$, and correspondingly,
$o(\ola)=k$, where $k\in\{1, 2, 3, 4\}$, and so $a^k\in M_X$.
Since $M=\lg a^i\rg \lg c\rg$ and $M_X=\lg a^i\rg \lg c^r\rg$, we know that $\lg a^i\rg =\lg a^k\rg$,
which implies that $i\in \{1, 2, 3, 4\}$.
Clearly,  if $\ox=\ZZ_2$, then $M_X=M$; if $\ox=D_8$ and $\o(\olc)=2$, then $M_X=\lg a^2\rg \lg c^2\rg$;
if $\ox =A_4$ and $\o(\olc)=3$, then $M_X=\lg a^2\rg \lg c^3\rg$;
if $\ox =S_4$ and $\o(\olc)=4$, then $M_X=\lg a^3\rg \lg c^4\rg $;
and  if $\ox=S_4$ and  $\o(\olc)=3$, then $M_X=\lg a^4\rg \lg c^3\rg$.
This is a contradiction.
Suppose that $\og$ is a semi dihedral group.
By the minimality of $\ox=\og\lg\olc\rg$,
with the same as the above, we also get a contradiction.
\vskip 3mm
{\it Step 2:  Find a contradiction.}
\vskip 3mm
Suppose that $M_X=1$.
Since both $\lg a\rg_X$ and $\lg c\rg_X$ are contained in $M_X$, we get $\lg a\rg _X=\lg c\rg_X=1$.
By Example~\ref{SD16or24}, we get $|G|\geq32$ as the minimality of $X$.
Now we shall show $G_X=1$.
Arguing by contradiction, assume that $G_X\ne1$.
If $|G_X|\gneqq 4$, then by $G=\lg a, b\rg \cong SD_{8n}$ we get $\lg a\rg _X\ne 1$, a contradiction.
So $|G_X|\le 4$.
Since $G_X\lhd G\cong SD_{8n}$, we know that $|G:G_X|\le 4$, which implies $|G|\le 16$,
contradicting to $|G|\geq32$.
Therefore, $G_X=1$.

Next, we consider the faithful (right multiplication) action of $X$
on the set of right cosets  $\O:=[X:\lg c\rg]$.
Since $X$ contains a regular subgroup $G\cong SD_{8n}\,(n\geq3)$,
by Lemma~\ref{SDI}, we get $X$ is imprimitive.
Pick a  maximal subgroup $H$ of $X$ which contains $\lg c\rg$ properly.
Then $H=H\cap X=(H\cap G)\lg c\rg=\lg a^s, b_1\rg \lg c\rg\lneqq X$,
for some $b_1\in G\setminus \lg a\rg$ and some $s$.
Note that the order of $b_1$ is either 2 or 4, and $\lg a^s\rg=H\cap\lg a\rg$.
Using  the same argument as that in Step 1, one has $a^s\in H_X$.  Set $\ox=X/H_X$.
Consider the faithful primitive action of $\ox$ on $\O_1:=[\ox:\oh]$,
with a cyclic regular subgroup of $\lg\ola\rg$, where $|\O_1|=s$.
By Proposition~\ref{Burnside}, a cyclic group of composed order is a Burnside group,
we know that either $s$ is a prime $p$ such that $\ox\le \AGL(1,p)$
or $s$ is composite such that $\ox$ is 2-transitive on $\O_1$.
In what follows, we consider these two cases, separately.
\vskip 3mm
Case (1):  $a^s=1$.
\vskip 3mm
In this case, we get that $b_1$ is an involution as $H\cap\lg a\rg=\lg a^s\rg=1$.
Replacing $b_1$ with $b$, we know  $H=\lg c\rg \rtimes \lg b\rg $ and $X=\lg c, b\rg .\lg a\rg $.
Since $\o(a)=4n$, we know that $s$ is composite.
Then $\ox $ is 2-transitive on $\O_1$.
By Proportion {\ref{cd}}, $\ox\le \AGL(l,q)$ for some prime $q$,
which contains a cyclic regular subgroup $\lg \ola\rg $ of order $q^l$.
By Lemma~\ref{AGL}, $\ox\cong S_4$ and $\o(\ola)=4$ so that $\o(a)=4$ (as $H_X\le \lg b, c\rg$),
which implies $|G|=8$, contradicting with $|G|\geq32$.
\vskip 3mm
Case (2): $a^s\neq1$.
\vskip 3mm
Firstly, show $s=p$, a prime.
To do that, recall $\ox=X/H_X$, $\oh=H/H_X$ and $\O_1:=[\ox:\oh]$.
Arguing by contradiction, assume that $s$ is composite.
Then $\ox$ is 2-transitive on $\O_1$, with a cyclic regular subgroup $\lg\ola\rg$.
Since $a^s\in H_X$, we get $H_X\ne 1$ and of course $H_X\nleqq \lg c\rg$.
Suppose that  $\lg a^j\rg \lg c\rg \le H$. Then $a^j\in M$.
Using the same arguments as that in  the first line of Step 1, we get $a^j\in M_X=1$.
Therefore, there exists an $l$ such that $bc^l\in H_X$, which implies $\oh=\lg \olc\rg$.
Then $\ox=\lg\ola\rg\lg\olc\rg$, a product of two cyclic subgroups, cannot be isomorphic to $S_4$.
But by Lemma~\ref{AGL}, we get $\ox\cong S_4$, a contradiction.
Therefore,  $s=p$, a prime and $X/H_X\le \AGL(1, p)$.

Since $a^{2n}\in H$, we can choose an involution in $H\setminus\lg a\rg$, namely $b$ for convenience, such that $H=\lg a^s,b\rg$.
Secondly, we consider the quotient group $\oh:=H/\lg c\rg_H=\overline{\lg c\rg \lg a^p,b\rg}$,
taking into account $s=p$, a prime.
Then $\lg \olc \rg_{\oh}=1$ and $o(\ola^p)=o(a^p)$.
Let $H_0/\lg c\rg _H=\lg \ola^{pj}\rg \lg \olc\rg $ be the biggest subgroup of $\oh$
containing $\lg \olc\rg$  and contained in the subset $\lg \ola^p\rg \lg \olc\rg$.
Since $|\oh|<|X|$ and $\og$ is either a dihedral group or a semi dihedral group, by Lemma~\ref{??} or the induction hypothesis on $\oh$,
we know that $H_0/\lg c\rg _H=\lg \ola^{pk}\rg \lg \olc\rg$, for one of $k$ in $\{1,2, 3, 4\}$,
which implies $\lg a^{pk}\rg\lg c\rg\lg c\rg_H=\lg c\rg\lg a^{pk}\rg\lg c\rg_H=\lg c\rg \lg a^{pk}\rg$,
giving $\lg a^{pk} \rg \lg c\rg \le H\le X$.
Since $M_X=1$, we have $a^{pk}=1$ where $ k\in \{1, 2, 3, 4\}.$
Note that $a^p\neq 1$ and $4\di\o(a)$.
Thus, we have that the order of $a$ is $4p$.
Therefore, only  the following  groups are remaining:
        $G=SD_{8p}$, where $p$ is a prime and $p\geq5$.
But by Lemma~\ref{n=p}, we get $G\neq SD_{8p}$, a contradiction.
\qed

\section{Proof of Theorem~\ref{main2}}
The proof of Theorem~\ref{main2} consists of the following two lemmas.

\begin{lem}\label{G_X}
Suppose that $\lg c\rg_X=1$ and $M=\lg a\rg \lg c\rg$.
Then $\lg a\rg _X\ne 1$.
\end{lem}
\demo
Since $\lg c\rg_X=1$,  by Proposition~\ref{cyclic}, we have $m< |G|$.
So $S:=G\cap G^{c}\ne 1$, otherwise $|X|\ge (8n)^2\gneqq |X|$.
Let $M=\lg a\rg \lg c\rg $, where $\o(a)=4n\geq8$ and $\o(c)=m$.
Arguing by contradiction, assume that $\lg a\rg _X=1.$
If $n\ge m$, then by Proposition~\ref{ab}, $\lg a\rg_M\ne 1$
and then  $\lg a\rg_X\ne 1$ is a contradiction.
So in what follows, we assume that  $n+1\le m$.

Since $\lg c\rg_M\ne 1$, we take  $z:=c^{\frac mp}\le \lg c\rg_M$ for a prime $p$.
Since $\lg c\rg_X=1$, we know that $\lg z^b\rg \ne \lg z\rg $ so that
$N:=\lg z\rg\times \lg z^b\rg \lhd X$.
Then $N$ contains an element $a^xc^y$ for some $x\neq0$ and $y$.
Set $a_1=a^{\frac {4n}p}$.

We claim that $z\in Z(M)$.
Suppose that $p=2$. Then $z\in Z(M)$, as desired.
Suppose that $p$ is odd. Then $p\di n$ as $\lg c\rg_X=1$.
Let $N\le P\in \Syl_p(M)$.
By Proposition~\ref{matecyclic}, $P$ is a  metacyclic group and
so we know that $N=\lg a_1\rg\times\lg z\rg $.
Since $\lg z\rg =\lg c\rg_M$, we may set $z^a=z^i$ and $z^b=a_1^jz^l$, where $j\ne 0$.
Then
$$(z^a)^b=(z^b)^{a^{2n-1}}=(a_1^jz^l)^{a^{2n-1}}=a_1^jz^{li^{2n-1}}\quad
{\rm and} (z^i)^b=a_1^{ji}z^{li},$$
which implies $a_1^{j(i-1)}=1$, that is $i=1$, and so  $z\in Z(M)$ again.

Since $z\in Z(M)\lhd X$, we get $z^b\in Z(M)$ which implies $N\leq Z(M)$.
Thus $a^xc^y\in Z(M)$, which implies $a^x\in Z(M)$, and then $\lg a^x\rg \lhd X$ for some $x\neq0$
is a contradiction.

\begin{lem}\label{SDa^2}
Suppose that $G=SD,\,X=X(SD),\,M=\lg a\rg \lg c\rg$ and $\lg c\rg_X=1$.
Then $\lg a^4\rg\lhd X$.
\end{lem}
\demo
Take a minimal counter-example $X$.
In the following Step 1, we show that the possible groups for $G$ are $SD_{8p^k}$,
where $p$ is a prime and $k\geq2$;
and in Step 2, we show that $G$ cannot be these groups.
\vskip 3mm
{\it Step 1: Show that the possible groups for $G$ are  $SD_{8p^k}$, where $p$ is an odd prime and $k\geq2$.}
\vskip 3mm
By Lemma~\ref{G_X}, we know $\lg a\rg_X\neq1$.
Suppose that $a^{2n}\in\lg a\rg_X$.
Then $\lg a^{2n}\rg\lhd X$ and $\lg a^{2n}\rg\lg c\rg\leq X$.
Set $\lg a^{2n}\rg\lg c_0\rg:=(\lg a^{2n}\rg\lg c\rg)_X$.
Since $\lg a^{2n}\rg\lg c_0\rg\lhd X$ and $\lg c\rg_X=1$, we get $c_0^2=1$.
Consider $\ox:=X/\lg a^{2n}\rg\lg c_0\rg=\og\lg \olc \rg$.
Noting that $\og$ is a dihedral group, $\lg\olc\rg_{\ox}=1$ and $\lg\ola\rg\lg\olc\rg\leq\ox$,
by Lemma~\ref{Da^2}, we get $\lg\ola^2\rg\lhd\ox$.
Then $\lg a^4\rg\char\lg a^2\rg\rtimes\lg c_0\rg\lhd X$.
Therefore, we get $\lg a^4\rg\lhd X$, as desired.
So in what follows, we assume that $a^{2n}\notin\lg a\rg_X$.
Then we know that the order of $\lg a\rg_X$ is odd and $\lg a\rg_X<\lg a^4\rg$ by minimality of $X$.

Let $p$ be the maximal prime divisor of $|\lg a\rg_X|$ and
set $a_0=a^{\frac {4n}p}\in\lg a\rg_X<\lg a^4\rg$.
Set $\ox=X/\lg a_0\rg=\og\lg \olc \rg$ and $\lg \olc \rg_{\ox}=\lg \olc_0\rg$.
Since $p$ is odd, we know $\og$ is a semi dihedral group.
(i) Suppose  that  that $\lg\olc\rg_{\ox}=1$.
Then by the minimality of $X$ we  get $\lg\ola^4\rg\lhd \ox$,
which implies $\lg a^4\rg \lhd X$, a contradiction.
(ii)Suppose  that  that $\lg \olc\rg \lhd \ox$.
Then $\ox/C_{\ox}(\lg \olc\rg)\le\Aut(\lg \olc\rg)$,
which is abelian and so $\ox'\le C_{\ox}(\lg \olc\rg)$.
Then $\ola^4\in \og'\leq \ox'\leq C_{\ox}(\lg \olc\rg)$,
that is $[a^4, c]\in\lg a_0\rg$, which implies $\lg a^4\rg \lhd X$, a contradiction.
By (i) and (ii), we have $1\ne \lg \olc\rg_{\ox}=\lg\olc_0\rg<\lg \olc\rg$.
Reset
$$K=\lg a_0\rg\rtimes\lg c_0\rg,\,\ox=X/K=\og\lg \olc\rg~{\rm and}~H=\lg a^4, c_0\rg.$$
If $\o(a_0)<\o(c_0)$, then $\{1\}\subsetneqq \lg c_0^j\rg=Z(K)\lhd X$, for some $j$, is a contradiction. Therefore, $1< \o(c_0)\le \o(a_0)$.
Then  we have the following two cases:
\vskip 3mm
{\it Case 1: $K=\lg a_0\rg \rtimes \lg c_0\rg \cong \ZZ_p\rtimes \ZZ_r$ is a Frobenius group, where $r\ge 2$.}
\vskip 3mm
In this case, by the minimality of $X$, we have $H/K=\lg \ola^4\rg\lhd \ox$,
that is $H=\lg a^4\rg \rtimes \lg c_0\rg \lhd X$.
Since $K\lhd X$, we know that $\lg a^4\rg/\lg a_0\rg$ and $\lg c_0\rg \lg a_0\rg/\lg a_0\rg$ are normal
in $H/\lg a_0\rg$.
Then $[a^4, c_0]\le \lg a_0\rg $.
So one can write
  $$H=\lg a^4,c_0|a^{4n}=c_0^r=1,\,(a^4)^{c_0}=a^4a_0^j\rg.$$
Let $P\in\Syl_p(H)$. Then $P\char H\lhd X$ so that $P\leq \lg a\rg_X$.
Clearly, one can check $Z(H)=\lg a^{4p}\rg$, which implies $\lg a^{4p}\rg\lhd X$.
Then $\lg a^{4p}\rg\leq \lg a\rg_X$.
Note that $\lg a^{4p},P\rg=\lg a^{4p}, a^{4n/p^k}\rg =\lg a^4\rg $, where $p^k\di\di 4n$,
so that  $a^4\in \lg a \rg_X$ is a contradiction again.
\vskip 3mm
{\it Case 2: $K=\lg a_0\rg \times \lg c_0\rg \cong \ZZ_p^2$.}
\vskip 3mm
With the same reason as that in Case 1, we have $H=\lg a^4\rg \rtimes \lg c_0\rg \lhd X$.

Let $H_1$ be the $p'$-Hall subgroup of $H$.
We get that $H_1$ is also the $p'$-Hall subgroup of $\lg a^4\rg$ as $\o(c_0)=p$.
Then $H_1\lhd X$, which implies $H_1\leq \lg a\rg_X$.
Suppose  that  $H_1\neq1$.
Let $a_1$ be an element of order $q$ in $H_1$, where $q<p$ is an odd prime as the maximality of $p$ and $|\lg a\rg|_X$ is odd.
Consider $\ox:=X/\lg a_1\rg=\og\lg c\rg $.
Similarly, we have $1\ne\lg\olc\rg_{\ox}:=\lg\olc_2\rg<\lg\olc\rg$ and
$H_0:=\lg a^4\rg\rtimes\lg c_2\rg\lhd X$.
Let $P\in \Syl_p(H_0)$.
Then $P\char H$ and so $P\lhd X$, which implies $P\leq \lg a\rg_X$.
Noting $\lg H_1,P\rg=\lg a^4\rg$, we therefore get $a^4\in \lg a\rg_X$, a contradiction.
So $H_1=1$, which means that $G$ is $SD_{8p^k}$ where $p$ is an odd prime and $k\geq2$.

\vskip 3mm
{\it Step 2: Show that the possible values of $m$ are $pq^e$, for a prime $q$ (may be equal to $p$)  and an integer $e$.}
\vskip 3mm
Arguing by contradiction, assume that $m=pq^em_1$ where $e\geq1$, $m_1\ne1$ and $q\nmid m_1$.
Recall $a_0=a^{4p^{k-1}}$, $c_0=c^{\frac mp}$ and set $a_2=a^4$.
Then $H=\lg a_2,c_0\rg=\lg a_2\rg\rtimes\lg c_0\rg.$
Note that $H$ is a $p$-group and $\lg a^{2p}\rg=\mho_1(H)\char H\lhd X$. Thus $\lg a^{4p}\rg\lhd X$.
Since $H\lhd X$, we get $H\lg c\rg\leq X$ and consider $\ox=X/{(H\lg c\rg)_X}=\og\lg\olc\rg$.
Note that $\og\cong\ZZ_4\times\ZZ_2$, $\lg\olc\rg_{\ox}=1$ and $\lg\ola\rg\lg\olc\rg\leq\ox$.
Then by Example~\ref{Z4Z2},
in the isomorphic sense, we get the form $\ox$ as shown in Table~\ref{table2}.

\begin{table}
  \center   \caption {The forms of $\ox$}\label{table2}
  \begin{tabular}{cccc}
  \hline
  Case & $\ox$\\
  \hline
  1&$\lg\ola,\olb,\olc\di\ola^4=\olb^2=\olc=1,\,\ola^{\olb}=\ola\rg$ \\
  2&$\lg \ola,\olb,\olc\di \ola^4=\olb^2=\olc^2=1,\ola^{\olb}=\ola^{\olc}=\ola,\olb^{\olc}=\ola^2\olb\rg$\\
  3&$\lg \ola,\olb,\olc\di \ola^4=\olb^2=\olc^2=1,\ola^{\olb}=\ola,\ola^{\olc}=\ola^3,\olb^{\olc}=\olb\rg$ \\
  4&$\lg\ola,\olb,\olc\di\ola^4=\olb^2=\olc^4=1,
   \ola^{\olb}=\ola^{\olc}=\ola,\olc^{\olb}=\ola\olc^3\rg$\\
  \hline
  \end{tabular}
  \end{table}
Then we have the following four cases:

(1) Suppose that $\ox=\lg\ola,\olb,\olc\di\ola^4=\olb^2=\olc=1,\,\ola^{\olb}=\ola\rg$.
Then in $\ox=X/H=\lg\olc\rg\rtimes\og$,
we get $\lg\olc^{q^e}\rg\lhd \ox$ and $\lg\olc^{m_1}\rg\lhd \ox$, which implies
$X_1=\lg a,b\rg\lg c^{q^e}\rg<X$ and $X_2=\lg a,b\rg\lg c^{m_1}\rg<X$.
By the minimality of $X$, we get $\lg a_2\rg\lhd X_1$ and $\lg a_2\rg\lhd X_2$.
Thus, $\lg a_2\rg\lhd\lg X_1,X_2\rg=X$ is a contradiction.
Therefore, $m=pq^e$. Suppose that $e=0$. Then $c_0=c$.
Consider $X_3=\lg a^2,b\rg\lg c\rg$.
If $\lg c\rg_X\neq1$, then $\lg c\rg\lhd X_3$ and we get $\lg a_2\rg\lhd X$, as $H=\lg a_2\rg\rtimes\lg c\rg$, a contradiction.
So $\lg c\rg_X=1$. Since $\lg a^2,b\rg$ is a dihedral group, by Lemma~\ref{Da^2}, we get $\lg a_2\rg\lhd X_3$, which implies $\lg a_2\rg\lhd X$, a contradiction.
Therefore, $m=pq^e$, where $q$ is a prime and $e\ge1$.

(2) Suppose that $\ox=\lg \ola,\olb,\olc\di \ola^4=\olb^2=\olc^2=1,\ola^{\olb}=\ola^{\olc}=\ola,\olb^{\olc}=\ola^2\olb\rg$.
Then $2\di m$ and we set $q=2$.
Consider $\ox=X/H$.
Then one can check $X_1=\lg a,b\rg\lg c^{2^e}\rg<X$ and $X_2=\lg a,b\rg\lg c^{m_1}\rg<X$.
By the minimality of $X$, we get $\lg a_2\rg\lhd X_1$ and $\lg a_2\rg\lhd X_2$.
Thus, $\lg a_2\rg\lhd\lg X_1,X_2\rg=X$ is a contradiction.

(3) Suppose that $\ox=\lg \ola,\olb,\olc\di \ola^4=\olb^2=\olc^2=1,\ola^{\olb}=\ola,\ola^{\olc}=\ola^3,\olb^{\olc}=\olb\rg$.
Then $2\di m$ and we set $q=2$.
Consider $\ox=X/H$.
Then one can check $X_1=\lg a,b\rg\lg c^{2^e}\rg<X$ and $X_2=\lg a,b\rg\lg c^{m_1}\rg<X$.
By the minimality of $X$, we get $\lg a_2\rg\lhd X_1$ and $\lg a_2\rg\lhd X_2$.
Thus, $\lg a_2\rg\lhd\lg X_1,X_2\rg=X$ is a contradiction.

(4) Suppose that $\ox=\lg\ola,\olb,\olc\di\ola^4=\olb^2=\olc^4=1,
   \ola^{\olb}=\ola^{\olc}=\ola,\olc^{\olb}=\ola\olc^3\rg$.
Then $4\di m$ and we set $q=2$.
Consider $\ox=X/H$.
Then one can check $X_1=\lg a,b\rg\lg c^{2^e}\rg<X$ and $X_2=\lg a,b\rg\lg c^{m_1}\rg<X$.
By the minimality of $X$, we get $\lg a_2\rg\lhd X_1$ and $\lg a_2\rg\lhd X_2$.
Thus, $\lg a_2\rg\lhd\lg X_1,X_2\rg=X$ is a contradiction.

Therefore, we get $m=pq^e$, $e\geq1$, and in particular, if $X/(H\lg c\rg)_X$ is the forms of Case 2,3,4 in Table~\ref{table2}, then $q=2$.
\vskip  3mm
{\it Step 3: Exclude the case $m=pq^e$, for a prime $q$  and an integer $e\geq1$.}
\vskip 3mm

Set $X_4:=(H.\lg c^q\rg).\lg a, b\rg =\lg a,b\rg\lg c^q\rg < X$.
By the minimality of $X_4$, we get $\lg a_2\rg\lhd X_4$,
that is $X_4=(\lg a_2\rg\rtimes\lg c^q\rg).\lg a, b\rg$ and $X_4/(\lg a_2\rg\rtimes\lg c^q\rg)$ is abelian.
Clearly, $\lg a_2\rg=G'\le X_4'\le \lg a_2, c^q\rg$.
So set $X_4'=\lg a^2, c_3\rg $ for some $c_3\in \lg c^q\rg$.
Since $\lg a_2^p\rg\lhd X$, by Proportion~\ref{NC},
both $X/C_{X}(\lg a_2^{p}\rg)$ and $X_4/C_{X_4}(\lg a_2\rg)$ are abelian,
which implies that $X'\leq C_X(\lg a_2^{p}\rg)$ and $X_4'\leq C_{X_4}(\lg a_2\rg)$.
Then $X_4'$ is abelian as $\lg a_2\rg\leq X_4'$.
The $p'$-Hall subgroup of $X_4'$ is normal, contradicting with $\lg c^q\rg_{X_4}=1$,
meaning that $X_4'$ is an abelian $p$-group.

Set $L:=H\rtimes\lg a^{2p^k}\rg=\lg a^2\rg\lg c_0\rg< X_4$.
We claim that $L\ntrianglelefteq X$. Arguing by contradiction, assume that $L\lhd X$.
If $H$ is abelian, then we get that either $\lg a_2\rg =Z(L)\char L\lhd X$, a contradiction;
or $L$ is abelian, forcing $\lg a^{2p^k}\rg\char L\lhd X$, a contradiction again.
Therefore, $H$ is non-abelian.
Note that $X_4'=\lg a_1, c_3\rg $ for $c_3\in \lg c^q\rg$.
If $c_3\ne1$, then  $c_0\in\lg c_3\rg\leq X_4'$ as $\o(c_0)=p$,
which implies that $H=\lg a^2,c_0\rg$ is abelian, a contradiction.
Therefore, $X_4'=\lg a^2\rg$, which implies $L=\lg a^2\rg\rtimes\lg c_0\rg$.
Noting that $\lg a^{2p^k}\rg\char L\lhd X$,
we get $\lg a^{2p^k}\rg\lhd X$, contradicting to $a^{2p^k}\notin\lg a\rg_X$.
\vskip 3mm
{\it Case (1): $q$ is an odd prime.}
\vskip 3mm
Suppose that $q$ is an odd prime. Then $X/(H\lg c\rg)_X$ is in Table~\ref{table2}.
In $\ox=X/H=(\lg\olc\rg\rtimes\lg \ola^2\rg).\lg\ola, \olb\rg$,
we get that $\olc^{\ola^2}=\olc^{-1}$ as $q$ is odd,
which implies $\lg a^2, c^{qp}\rg\leq X_4'\leq \lg a^2\rg\lg c^q\rg$.
Note that  $X_4'$ is the abelian $p$-group.
Thus either $q\ne p$ and $e=1$; or $q=p$.
Suppose that $q\ne p$ and $e=1$, that is $\o(c)=pq$.
Consider $M=\lg a\rg\lg c\rg\lhd X$.
Then by Proportion~\ref{mateabel}, $M'$ is abelian.
Note that $\lg c\rg_X=1$ and $\lg a\rg_X$ is the $p-$group.
Thus $M'$ is an abelian $p$-group with the same argument as the case of $X_4'$.
Noting $\lg a_2\rg\lg c^p\rg$ is the $p'$-Hall subgroup of $M$,
we get $[a_2,c^p]\in \lg a_2\rg\lg c^p\rg\cap M'=1$, which implies $\olc^{\overline{a_2}}=\olc$ in $\ox=X/H$, a contradiction.
So in what follows, we assume $q=p$, that is $\o(c)=p^{e+1}$.
Note that $\olc^{\ola^2}=\olc^{-1}$ in $\ox=X/H$ and $\lg a_2\rg\leq X_4'$.
Thus $X_4'$ is either $\lg a_2\rg\lg c^p\rg$ or $\lg a_2\rg$,
noting $X_4'=\lg a_2\rg$ only happens when $e=1$.

Suppose that $X_4'=\lg a_2\rg\lg c^p\rg$. Note that $H\leq X_4'$.
Thus $H=\lg a_2\rg\rtimes\lg c_0\rg$ is abelian.
Note that both $\lg a_2\rg$ and $\lg c\rg$ are $p$-groups and $X=(H.\lg c\rg)\rtimes\lg a^{p^k}, b\rg$.
Set $$a_2^{c}=a_2^{s}c_0^{t}\quad{\rm and}\quad c^b=a_2^{u}c^v,$$
where $s\equiv1(\mod p)$ and $p\nmid v$.
Then for an integer $w$, we get
$$a_2^{c^w}=a_2^{x_1}c_0^{wt}\quad {\rm and}\quad c_0^b=a^{x_2}c_0^{v}$$
for some integers $x_1$ and $x_2$.
Since $(a_2^{c})^b=(a_2^{s}c_0^{t})^b$, there exist some integers $x$ and $y$ such that
$$(a_2^{c})^b=(a_2^{-1})^{c^v}=a^{x}c_0^{-vt}\quad{\rm {and}}\quad
(a_2^sc_0^{t})^b=a^yc_0^{vt},$$
which gives $t\equiv0(\mod p)$.
Then $\lg a_2\rg\lhd X$ is a contradiction.

Suppose that $X_4'=\lg a^2\rg$. Then $\o(c)=p^2$ and $c_0=c^p$.
Note that $X_4=G\rtimes\lg c^q\rg=H\rtimes\lg b\rg$.
Thus $[c_0,a^{2p^k}]=1$.
Set $c^{a^{2p^k}}=a^xc^{-1+yp}$ as $\olc^{\ola^2}=\olc^{-1}$ in $\ox=X/H$.
Then
$$c=c^{(a^{2p^k})^2}=a^x(a^xc^{-1+yp})^{-1+yq}=a^xc^{1-yp}a^{-x}c^{yp},$$
which implies $(a^x)^{c^{1-yp}}=a^x$.
Then $[a^x,c]=1$.
Note that $c_0=c^p$ and $c_0=c_0^{a_2}=(c^{a^{2p^k}})^p=a^{x}c^{-p}$ for some $x$.
Thus we get $c_0^2=1$, contradicting with $\o(c_0)=p$.

\vskip 3mm
{\it Case (2): $q=2$.}
\vskip 3mm
In this case, consider $\ox=X/H=(\lg\olc\rg\rtimes\lg \ola^2\rg).\lg\ola, \olb\rg$.
If $e=1$, that is $\ox=(\lg\olc\rg\times\lg \ola^2\rg).\lg\ola, \olb\rg$, then we get
$L=\lg a^2\rg\lg c_0\rg\lhd X$, a contradiction.
So $e\geq2$.
Then we get that $X_4=\lg a,b\rg\lg c^2\rg\lhd X$ and $\olc^{\ola^2}$ is either $\olc^{-1}$ or $\olc^{\pm1+2^{e-1}}$ as $L\ntriangleleft X$.
Note that $a_2\in X_4'$, $X_4'\char X_4\lhd X$ and $X_4'$ is the abelian $p$-group.
Thus $X_4'=H$, which implies that $H$ is ablian.

We shall show $\olc^{\ola_2}=\olc^{1+2^{e-1}}$ in $X/H$.
Arguing by contradiction, assume that  $\olc^{\ola^2}$ is either $\olc^{-1}$ or $\olc^{-1+2^{e-1}}$.
Then $(\olc^2)^{\ola_2}=\olc^{-2}$, which implies $\lg a_2,c^4\rg\leq X_4'$.
Then $c^4\in\lg c_0\rg$, which implies $e=2$ as $e\geq2$.
But when $e=2$, we get $\olc^{\ola_2}=\olc^{-1}=\olc^{1+2}$, a contradiction.
Therefore, $\olc^{\ola^2}=\olc^{1+2^{e-1}}$ in $X/H$.
By Proportion~\ref{sylowp}, we get that $\lg a^{p^k}, a^ib\rg\lg c^p\rg$ and $\lg a^{p^k}\rg\lg c^p\rg$ are  Sylow 2-subgroups of $X$ and $M$.
Note that $\lg a^{p^k}, a^ib\rg\lg c^p\rg=(\lg c^p\rg\rtimes \lg a^{2p^k}\rg).\lg a, a^ib\rg$.
Then $(c^p)^{a^{2p^k}}=c^{p+2^{e-1}p}$ and $e\ge2$, which implies $c^{2^{e-1}p}\in M'$.

Suppose that $e=2$. Then $(c^p)^{a^{2p^k}}=c^{-p}$.
Since $|\lg a^{p^k}\rg|\geq|\lg c^p\rg|$, by Proportion~\ref{ab},
we get $\lg a^{2p^k}\rg\lhd\lg a^{p^k}\rg\lg c^p\rg$.
Since $a^{2p^k}$ is an involution, we get $[a_2,c^p]=1$, a contradiction.
So in what follows, we assume $e>2$. And set $a_3=a^{2p^k}$.

Noting $\lg c^p\rg\rtimes\lg a_3\rg=\lg a_3,c^p|a_3^2=c^{2^ep}=1,\,(c^p)^{a_3}=c^{(1+2^{e-1})p}\rg$,
there are only three involutions in $\lg c^p\rg\rtimes\lg a_3\rg$: $a_3,c^{2^{e-1}p}$ and $a_3c^{2^{e-1}p}$.
Since $\lg a^{p^k}, a^ib\rg\cong\ZZ_4\times\ZZ_2$, by Table~\ref{table2}, we get
$\lg c^{2p}\rg\lhd\lg a^{p^k}, a^ib\rg\lg c^p\rg$ for Case 1,2,3.
Then we get $[a^{p^k}, c^{2^{e-1}p}]=[a^ib, c^{2^{e-1}p}]=1$
Recall $M=\lg a\rg\lg c\rg$. By Proportion~\ref{mateabel}, $M'$ is abelian.
Let $M_2$ be the Sylow 2-subgroup of $M'$.
Note that $M_2\char M'\char M\lhd X$, $c^{2^{e-1}p}\in M'$, $\lg c\rg_X=1$ and $a_3$ is an involution.
Thus we get $M_2\cong \ZZ_2^2$, which implies $M_2=\lg c^{2^{e-1}p}, a_3\rg$.
Consider $HM_2\leq X$.
Since $M_2\lhd X$, $H\lhd X$, $H\cap M_2=1$ and $p$ is odd prime,  we get $HM_2=H\times M_2\lhd X$,
which implies $a^2$ normalises $\lg c^{2^{e-1}p}\rg$.
Since  $X=\lg a,b,c\rg$ and $[a^{p^k}, c^{2^{e-1}p}]=[a^ib, c^{2^{e-1}p}]=1$, we get $\lg c^{2^{e-1}p}\rg\lhd X$, a contradiction.
So in what follows, we assume that Case 4 in Table~\ref{table2},
that is $X/(H\lg c\rg)_X=\lg \ola,\olb,\olc\di \ola^4=\olb^2=\olc^4=1,\ola^{\olb}=\ola^{\olc}=\ola,\olc^{\olb}=\ola\olc^3\rg.$

In this case, we know $e\geq2$ and set
\begin{eqnarray}\label{f4.1}
a^c=a^{1+4i}c^{4j},\,c^b=a^{1+4k}c^{3+4l}.
\end{eqnarray}
Noting that $X_1=\lg a,b,c^2\rg=\lg a,b\rg\lg c^2\rg<X$, by the minimality of $X$, we get $\lg a_1\rg\lhd X_1$.

Consider $\ox=X/\lg a_0\rg$.
Since $\lg\olc_0\rg\lhd\ox$ and $\oh=\lg \ola_1\rg\times\lg\olc_0\rg$,
we get $\lg\ola_1,\olc\rg\leq C_{\ox}(\lg \olc_0\rg)\lhd\ox$.
By Eq(\ref{f4.1}), we get $\lg\ola,\olc\rg\leq C_{\ox}(\lg \olc_0\rg)$.
Then we know $H\lg a\rg=\lg a\rg\rtimes\lg c_0\rg$.
Consider $\ox=X/H$.
Since $\lg\olc^4\rg\lhd\ox$,
we get $\lg\ola,\olc\rg\leq C_{\ox}(\lg \olc^4\rg)$ with the same as the above.
By Eq(\ref{f4.1}), we get $\ola^{\olc}=\ola\olc^{4j}$.
Since $\o(\ola)=4$, we know that $\o(\olc^{4j})$ is 1, 2 or 4.

Arguing by contradiction, assume that $\o(\olc^{4j})\leq2$.
Then we get $(\ola^2)^{\olc}=\ola^2$, which implies $H\lg a^2\rg\lhd X$.
Since $H\lg a^2\rg=\lg a^2\rg\rtimes\lg c_0\rg$, we get
$$\lg a^{p^k}\rg\char H\lg a^2\rg\lhd X,$$
which implies $a^{p^k}\in \lg a\rg_X$.
Since the order of  $\lg a\rg_X$ is odd, we get a contradiction.
Thus, $\o(\olc^{4j})=4$, which implies $16\di m$ and $m_1\geq4$.
Set $c_2=c^{2^{m_1-2}}$ and $c_3=c^{2^{m_1-2}p}$.
Then we can reset $a^c=a^{1+4i}c_2^{j}$.
Moreover, we get $(\ola^2)^{\olc^2}=\ola^2$ in $\ox=X/H$.

Recall $K=\lg a_0\rg\times\lg c_0\rg\cong\ZZ_p^2$.
Consider $\ox=X/K=\og\lg\olc\rg$.
Since $\lg\olc\rg_{\ox}=1$ and $\og$ is a semi dihedral group,
by the minimality of $X$, we get $\lg \ola_1\rg\lhd\ox$.
Note that $\lg\ola\rg\leq C_{\ox}(\lg \ola_1\rg)\lhd \ox$.
Since $\lg\ola_1\rg$ is cyclic, we know $\ox/C_{\ox}(\lg \ola_1\rg)$ is abelian.
By Eq(\ref{f4.1}), we know that $\olc^2\in C_{\ox}(\lg \ola_1\rg)$.
By Proportion~\ref{sylowp}, we know that $\lg a^{p^k},b\rg\lg c^p\rg$ is a  Sylow 2-subgroup of $X$.
By Eq(\ref{f4.1}) again, we can set $(a^{p^k})^{c^p}=a^{p^k}c_3^j.$
Since $(\ola^2)^{\olc^2}=\ola^2$ in $\ox=X/H$, we get $(a^{2p^k})^{c^{2p}}=a^{2p^k}.$
Then we know $\lg a^{2p^k}\rg\lhd X_1$.
Since $a^{2p^k}$ is an involution, we know $a^{2p^k}\in Z(X_1)$.
We claim that $a^{2p^k}$ is a unique involution in $Z(X_1)$.
Arguing by contradiction, assume that $a^xb^yc^z\in Z(X_1)$ is an involution.
Since $1=[a^xb^yc^z,c]=[a^xb^y,c]$, we get $(a^xb^y)^2=(c^z)^2=1$ as $(a^xb^yc^z)^2=1$.
If $y=0$, then $a^xb^yc^z=a^{2p^k}c_3^2\in Z(X_1)$, recalling $c_3^2=c^{2^{m_1-1}p}$, and
we get $c_3^2\in Z(X_1)$, as $a^{2p^k}\in Z(X_1)$, contradicting to $\lg c^2\rg_{X_1}=1$.
So in what follows, we assume that $y=1$.
Then $a^xb^yc^z=a^xbc_3^2\in Z(X_1)$, which implies $\overline{a^xbc_3^2}\in Z(X_1/K)$.
But $\ola_1^{\overline{a^xbc_3^2}}=\ola_1$ implies $\ola_1^{\olc_3^2}=\ola_1^{-1}$,
contradicting with $\olc^2\in C_{\ox}(\lg \ola_1\rg)$.
Thus, $a^{2p^k}$ is a unique involution in $Z(X_1)$.
Then $$\lg a^{2p^k}\rg\char X_1\lhd X,$$
which implies $a^{p^k}\in \lg a\rg_X$.
Since the order of  $\lg a\rg_X$ is odd, we get a contradiction again.
\qed
\vskip 3mm
Note that for $G\in\{D,Q\}$, if $\lg a\rg\lhd X(G)$ and $\lg c\rg_{X(G)}=1$, then $G\lhd X$.
We shall consider the case, that is $G=SD$, $\lg a\rg\lhd X(G)$ and $\lg c\rg_{X(G)}=1$,
and find something interesting that differs from the case of $G\in\{D,Q\}$.
Then we have the following lemma.
\begin{lem}\label{a}
Suppose that $G=SD,\,X=X(G),\,\lg c\rg _X=1$ and $M=\lg a\rg\lg c\rg$.
If $\lg a\rg\lhd X$, then $c^2\in N_X(G)$, and especially $G\lhd X$
if $b\in N_X(\lg a^2\rg\lg c\rg)$.
\end{lem}
\demo
$X=(\lg a\rg\rtimes\lg c\rg).\lg b\rg,$ and so we may write $a^c=a^i$ and $c^b=a^kc^j$.
If $j=1$, then $G\lhd X$. So in what follows, we assume that $j\ne 1$.
Since $b^2=1$, we get $c=c^{b^2}$.
Then
$$c=(c^b)^b=(a^kc^j)^b=a^{-k+2kn}(a^kc^j)^j=a^{2kn}c^j(a^kc^j)^{j-1}.$$

We claim that $k$ is odd. Arguing by contradiction, assume that $k$ is even.
Then $c^{1-j}=(a^kc^j)^{j-1}=(c^{j-1})^b,$ so that $b$ normalizes $\lg c^{1-j}\rg$.
Since $\ox=X/C_X(\lg a\rg)\le \Aut(\lg a\rg )$ which is abelian,
we get $\olc =\olc^{\olb}=\olc^j$, that is $c^{1-j}\le C_X(\lg a\rg)$ so that $[c^{1-j}, a]=1$.
Thus we get  $\lg c^{1-j}\rg \lhd X$.
It follows from $\lg c\rg _X=1$ that $j=1$, a contradiction.
Therefore, $k$ is odd.

Noting that $k$ is odd, we get $c=(c^b)^b=a^{2n}c^j(a^kc^j)^{j-1}.$
Since $a^{2n}$ is an involution of $X$ and $\lg a\rg\lhd X$, we get $a^{2n}\in Z(X)$.
Since $(c^2)^b=a^kc^ja^kc^j=a^{k(1+i^{-j})}c^{2j}$, we get that
$$c^2=((c^2)^b)^b=(a^{k(1+i^{-j})}c^{2j})^b=c^{2j}(a^{k(1+i^{-j})}c^{2j})^{j-1}.$$
With the same as the above, we get $2(j-1)\equiv0(m)$.
Then $(c^2)^b=a^{k(1+i^{-j})}c^{2}$, which implies that $c^2\in N_X(G)$ as $\lg a\rg\lhd X$, as desired.
\qed

\section{Conjecture}
Let $X(G)=G\lg c\rg$ be a group, where
$$G=\lg a,b\di a^n=1, b^2=a^t, a^b=a^r, r^2\equiv1(\mod n), t(r-1)\equiv0(\mod n)\rg\cong\ZZ_n.\ZZ_2$$
and $C$ is a cyclic group such that $G\cap C=1$.
Then $X=X(G)=G\lg c\rg=\lg a, b\rg \lg c\rg$.
Let $X$ contain a subgroup $M$ of the biggest order such that
$\lg c\rg\le M\subseteqq\lg a\rg\lg c\rg$.
Then we have the following conjecture.
\begin{conj}
Let $G=\lg a,b\di a^n=1, b^2=a^t, a^b=a^r, r^2\equiv1(\mod n), t(r-1)\equiv0(\mod n)\rg$ and $X=X(G)=G\lg c\rg$, where $\o(c)=m\ge 2$ and $G\cap \lg c\rg=1$.
Let $M$ be the subgroup of the biggest order in $X$ such that
$\lg c\rg \le M\subseteqq \lg a\rg \lg c\rg$.  Then one of items in Tables \ref{table3} holds.
\begin{table}
  \center \caption {The forms of $M$, $M_X$ and $X/M_X$}\label{table3}
  \begin{tabular}{cccc}
  \hline
  Case &$M$ & $M_X$  & $X/M_X$\\
  \hline
   1 & $\lg a\rg\lg c\rg$   & $\lg a\rg\lg c\rg$     &   $\ZZ_2$ \\
   2 &$\lg a^2\rg \lg c\rg$ & $\lg a^2\rg\lg c^2\rg$ &   $D_8$   \\
   3 &$\lg a^2\rg \lg c\rg$ & $\lg a^2\rg\lg c^3\rg$ &   $A_4$   \\
   4 &$\lg a^4\rg \lg c\rg$ & $\lg a^4\rg\lg c^3\rg$ &   $S_4$   \\
   5 &$\lg a^4\rg \lg c\rg$ & $\lg a^4\rg\lg c^4\rg$ &   $(\ZZ_4\times\ZZ_2)\cdot\ZZ_4$   \\
   6 &$\lg a^3\rg \lg c\rg$ & $\lg a^3\rg\lg c^4\rg$ &   $S_4$   \\
   \hline
  \end{tabular}
  \end{table}
\end{conj}

\end{document}